\documentclass[12pt,a4paper]{amsart}
\newtheorem{thm}{Theorem}[section]

\newtheorem{rem}[thm]{Remark}
\def\trace{\operatorname{\mathrm trace}\nolimits}
\def\Gal{\operatorname{\mathrm Gal}\nolimits}
\def\SL{\operatorname{\mathrm SL}\nolimits}
\def\SU{\operatorname{\mathrm SU}\nolimits}
\def\Sp{\operatorname{\mathrm Sp}\nolimits}
\def\Lie{\operatorname{\mathrm Lie}\nolimits}
\def\Frob{\operatorname{\mathrm Frob}\nolimits}
\title[ $\ell$-adic Poisson Formula and Endoscopy  for $p$-Adic Groups]{ $\ell$-adic Poisson Formula and Endoscopy  for $p$-Adic Reductive Groups}
\author{Do Ngoc Diep}
\address{Institute of Mathematics, Vietnam Academy of Sciences and Technology, 18 Hoang Quoc Viet Road, 10307 Hanoi, Vietnam}
\email{dndiep@math.ac.vn}
\date{\bf\today}

\begin{document}
\maketitle
\begin{abstract} For two distinguished prime $\ell$ and $p$,
 we prove a $\ell$-adic version of the Poisson formula for reductive $p$-adic groups. In order to do this we write an identity for the trace of regular representation and orbital integrals. Next we reduce them to orbital integrals for endoscopy groups and look at this  as the special value of $L$-function at $s=0$.
And finally show that it is equal to  the special value of motivic $L$-function at $s=0$.


\end{abstract}
\section{Introduction}  
The classical Poisson summation formula states that for any smooth compactly supported function $f$ on the real line $\mathbb R$, the following identity holds 
$$\sum_{n=-\infty}^{+\infty} f(x -n) = 2\pi\sum_{m=-\infty}^{+\infty}\hat{f}(m) e^{ imx},$$
where $\hat{f}(m) := \frac{1}{2\pi}\int_{-\pi}^{+\pi} f(x) e^{- i mx}dx$ is the Fourier transform of $f$. 

The right hand side is the inverse Fourier transform, which is uniformly convergent on $\mathbb R$. Writing out the finite sum
$$S_n(f)(x) =  2\pi\sum_{m=-n}^{+n}\hat{f}(m) e^{ imx},$$
we can interchange the order of summation and integration to have
$$S_n(f)(x) = 2\pi \sum_{m=-n}^{+n}\frac{1}{2\pi} \int_{-\pi}^{+\pi} f(x) e^{- i my}e^{imx}dy$$
$$=\sum_{m=-n}^{+n} \int_{-\pi}^{+\pi} f(x) e^{ i m(x-y)}dy$$
$$=\sum_{m=-n}^{+n} \int_{-\pi+2\pi m}^{+\pi+2\pi m} f(x-2\pi m) e^{ i m(x-y)}dy
= \sum_{n=-\infty}^{+\infty}f(x-n).$$
On distribution language  it sounds as
$$\sum_{n=-\infty}^{+\infty} \delta_n(x) = 2\pi \sum_{m=\infty}^{+\infty} e^{imx}$$
The Poisson summation formula plays an important role in mathematical physics and many other branches of mathematics where one manipulate the Fourier transforms.

Following the ideas  of J.-P. Labesse, in works \cite{diepquynh1}-\cite{diepquynh4} we show a version of Poisson formula as equality of the sum of traces of cuspidal automorphic representations of groups of split rank 1: $\SL(2,\mathbb R)\cong \Sp(2,\mathbb R) \cong \SU(1,1)$, and groups of split rank 2: $\SL(3,\mathbb R)$, $\Sp(4,\mathbb R)$ and $\SU(2,1)$, and  the corresponding sum of orbital integrals. 
In this paper we use some advanced results of J. Arthur, J.-L. Waldspurger, M. Harris and R. Taylor and others to prove the general case of reductive $p$-groups.

We use the notations from \cite{arthur},\cite{milne},\cite{harris}, \cite{colmez}.
In this paper we prove the following theorem:
\begin{thm}[The Main Result] \label{mainthm}
Let $G$ be a reductive group over number field $F$, $\Gal(F/\mathbb Q)$ the Galois group of finite  extension $F/\mathbb Q$.
For any $f$ from the Hecke algebra $\mathcal H(G)$ of bi-invariant functions with convolution,
$$\sum_{\pi\in \mathcal R} m(\pi)\trace(\pi(f)) =\sum a(\gamma) \mathcal O(\gamma, \hat f) = \sum_{\chi_{\ell,\pi}}a(\pi_{\ell,\pi}) L(0, \chi_{\ell,\pi},f),$$
where $\pi$ is running over the set $\mathcal R$ of all cuspidal automorphic representations of $G$, $m(\pi)$ is the multiplicity of representation $\pi$ in the regular representation, $a(\gamma)$ some coefficients, $\mathcal O(\gamma,\hat f)$ is the orbital integral of the Fourier transform of $f$,  $\pi(f)$ is the Hecke operator operating in representation $\pi$,
$\chi_{\ell,\pi}$ is running over the Haris-Taylor's Langlands parametrization of cuspidal automorphic representations of Galois group, $L(0, \chi_{\ell,\pi},f)$ are the special value  at $0$ of motivic L-functions $L(s, \chi_{\ell,\pi},f)$.
\end{thm}
As a straight corrolary we have
\begin{thm}
Let $G$ be a reductive group over number field $F$, $\Gal(F/\mathbb Q)$ the Galois group of finite  extension $F/\mathbb Q$.
For any $f$ from the Hecke algebra $\mathcal H(G)$ of bi-invariant functions with convolution,
$$\sum_{\pi\in \mathcal R} m(\pi)\trace(\pi(f)) =\sum a(\gamma) \mathcal O(\gamma, \hat f) $$ $$= \sum_{\chi_{\ell,\pi}}a(\pi_{\ell,\pi}) \det\left(1 - \chi_{\ell,\Pi}(\Frob_p)(f)\right),$$ where $Frob_p$ is the Frobenius element of Galois groups acting on the character space.
\end{thm}

The idea of the proof is as follows. First we use the J. Arthur's and Waldspurger's trace formula sum to reduce the sum of trace of subrepresentations of the regular representation to a sum of orbital integrals. Use the transfer process to reduce the problem to the  endoscopic version. At endoscopic level, the representations are multiple of one dimensional representatipons. For one dimensional representations, i.e. characters, there is a passage from these characters to the characters of maximal ablelian extension of field, i.e. Dirichlet characters. We reduce the interested sum to sum of Dirichlet L-functions. And finally we use the motivic interpolation for Dirichlet L-functions to have a sum of motivic L-funtions.
The particular cases of real ranks 1 and 2 were treated in \cite{diepquynh1}-\cite{diepquynh4}.

Let us describe the structure of the paper. In the next section \S2 we prove the main result: we describe the Arthur's theory of trace formula. Then in  we describe th Waldspurger's theory of local trace formula. Finally in we prove the stated main result. In \S3 we make some remark about the particular cases of the Poisson trace formula.

\section{Proof of the main theorem}
\subsection{Langlands' Parametrization of cuspidal automorphic representations}
The case of reductive groups over the real field is well-kown. The irreductive cuspidal automorphic representations are described by the Langlands parametrization. 

Let us remind the construction of cuspidal automorphic representations of $p$-adic reductive groups. Let $K$ be a finite extension of $\ell$-adic fields $\mathbb Q_\ell$ number field, $G = \underbar{G}(K)$ be the $K$-points of a split reductive group $\underbar{G}$ over $\mathbb Q$. Let $\mathfrak g = \Lie G$; more precisely,  for any place $v$, let $\mathfrak g_v = \Lie(G(K_v))$ be the Lie algebra $G(K_v)$. For any archimedian place $v$, let $U(\mathfrak g_v)$ be the universal envelloping algebra of $\mathfrak g_v$ and in the case of archimedian place $v$,  $\mathcal Z(\mathfrak g_v)$ - its centre which is isomorphic to the symmetric polynomial algebra  of the Cartan subalgebra $t_v\subset \mathfrak g_v$ , and $U_v\subset G(K_v)$ a maximal compact open subgroup. 
Following the well-known result of Harish-Chandra, there is one-to-one correspondence between the maximal ideals of $\mathcal Z(\mathfrak g_v)$ and the infinitesimal characters $\lambda_{\Pi_v}$ of algebraic smooth representations $\Pi_v$ of $G(K_v)$, which are irreducible $(\mathfrak g_v, U_v)$-modules with central character $\lambda_{\Pi_v}$. Let us consider the fundamental representation $\tau$ of $U_v$ in an $\ell$-adic vector space $V$. Because $\mathbb C \cong \bar{\mathbb Q}_\ell $ following the H. Weyl theorem, every irreducible representation of $U_v$ is an irreducible constituent of the symmetric $n$ power $S(\check{\tau})\otimes{\sqrt{\Delta}}$ of the contragradient (to the fundamental) representation $\check{\tau}$.

\subsection{Trace of $\ell$-adic Regular Representation}
Consider the case of real field $F=\mathbb R$:

The classical Poisson summation formula is the starting point for the Poisson formula for general reductive groups $G=G(\mathbb R)$ over real field.

Following  J. Arthur (\cite{arthur}, (29.7)),
for any $f\in \mathcal H(G)$ the trace $\trace(\mathcal R(f))$ can be computed as the sum:
$$\sum_{t\geq 0}\sum_{M\in \mathcal L} \frac{|W_M|}{|W|}\int_{\Phi_t(M,V,\zeta)}a^M(\pi) I_M(\pi,f)d\pi$$ $$=\sum_{M\in \mathcal L} \frac{|W_M|}{|W|}\sum_{\Phi_t(M,V,\zeta)}\sum a^M(\gamma)I_M(\gamma,f),$$ the coefficients $a(.)$  are defined in (\cite{arthur}, Theorem 29.4).

\subsection{Reduction to Orbital Integrals of Endoscopic Groups}
Againce in the case of the real field $\mathbb R$,
the idea is to reduce the intergrals to some integrals over endoscopic subgroups. This process is called the transfer:

J. Arthur in (\cite{arthur}, Corollary 29.9) showed that for any $f\in \mathcal H(G)$
$$\sum_{t\geq 0}\sum_{M\in \mathcal L} \frac{|W_M|}{|W|}\int_{\Phi_t(M,V,\zeta)}b^M(\phi)S_M(\Phi,f)d\Phi$$ $$=\sum_{M\in \mathcal L} \frac{|W_M|}{|W|}\sum_{\Phi_t(M,V,\zeta)}\sum b^M(\delta)S_M(\delta,f),$$ the coefficients $a(.)$ and $b(.)$ are defined in (\cite{arthur}, Theorem 29.4).

For $p$-adic reductive groups $G=G(\mathbb Q_p)$ some similar formula can be deduced from the work of J.-L. Waldspurger of twisted  local trace formula.

Following R. E. Kottwitz (\cite{kottwitz}, Theorem 10.4)
$\trace(\mathcal R(f))= J(f)= $ $$= \sum_{M\in\mathcal L} \frac{|W_M|}{|W|} \sum_{T\in \mathcal T_M}\frac{1}{n^M_T}\int_{\mathfrak t_{reg}}|D(X)| \int_{A_M\backslash G} f(g^{-1}Xg)\tilde{\mathcal V}_M(1,g; D)d\dot gdX$$

R. E. Kottwitz has showed in \cite{kottwitz} that the 
$$J(f_1,f_2) = J(\hat{f}_1, \check{f}_2),$$ where $\hat{f}$ is the Fourier transform and $\check{f}$ is the inverse Fourier transform with respect to an additive character.
From this, it is easy to see that the same result is true for the Fourier transform
$$J(f,1) = J(\hat{f}, \check{1}).$$

The endoscopic groups are the connected component of the identity of the centralizer of regular semisimple elements in the maximal tori. So the representations we are interested in are multiple of one-dimensional representations of the endoscopic subgroups $H$. Every characters factorizes through the abelianized group $H^{ab} = H/(H,H)$. The characters are irreducible and therefore it is factorized through the characters of multiplicative groups. 

\subsection{Reduction to Motivic $L$-Function}

Let us consider the complex field $\mathbb C\cong \bar{\mathbb Q}_\ell$ and consider the $\ell$-adic representations of our $p$-adic groups $G(\mathbb Q_p)$ in $\mathbb C\cong \bar{\mathbb Q}_\ell$. Certainly the Dirichlet $L$-series is highly divergent in the $\ell$-adic topology. 
Let us consider the characters of abelian endoscopic groups. Following the Class-Field Theory every such character factorizes through the Abelian Galois groups. These characters are corresponding to the Dirichlet characters of the maximal abelian extension $F^{ab}$. For Dirichlet character, we write out the motivic interpolation as motivic $L$-functions. Remark that the trace of Frobenius element can be consider as a distribution on the Hecke algebra. 
The following result of Artin is well-known,
see \cite{kowalski}
\begin{quote}
Let $K$ be a number field, $E/K$ a finite abelian extension with
Galois group $G$, let $\rho : G \to \mathbb C^\times$ be a Galois character and $L(\rho, s)$ the associated $L$-function.
Then there exists a unique primitive Hecke character $\chi$ of $K$, of modulus $\mathfrak m$ say, such that 
$$L(\rho, s) = L(\chi, s).$$
Actually, this identity holds locally, i.e. the $p$-component of the Euler products are equal
for all prime ideal $\mathfrak p$ of $K$, namely
$$\rho(\sigma_{\mathfrak p} ) = \chi(\mathfrak p)$$ for all $\mathfrak p$  coprime to $\mathfrak m$.
\end{quote}

Following  P. Colmez (\cite{colmez}, \cite{harris})
$$\trace(\pi(f))=L(0, \chi_{\ell,\pi},f).$$

The proof of our main Theorem \ref{mainthm} is therefore achieved. $\hfill\Box$

\section{Corollaries}

\begin{rem} The particular cases of the real field and low rank reductive groups of real rank 1 or 2 were treated in \cite{diepquynh1}-\cite{diepquynh4}.
\end{rem}

Following  P. Colmez (\cite{colmez}, \cite{harris}) for the automorphic Galois characters $\pi= \rho$ and the corresponding Hecke character $\chi_{\ell,\pi}$, 
$$\trace(\pi(f))=L(0, \chi_{\ell,\pi},f)=\det\left(1 - \chi_{\ell,\Pi}(\Frob_p)(f)\right).$$ We have the result as follows.

\begin{thm}
Let $G$ be a reductive group over number field $F$, $\Gal(F/\mathbb Q)$ the Galois group of finite  extension $F/\mathbb Q$.
For any $f$ from the Hecke algebra $\mathcal H(G)$ of bi-invariant functions with convolution,
$$\sum_{\pi\in \mathcal R} m(\pi)\trace(\pi(f)) =\sum a(\gamma) \mathcal O(\gamma, \hat f) $$ $$= \sum_{\chi_{\ell,\pi}}a(\pi_{\ell,\pi}) \det\left(1 - \chi_{\ell,\Pi}(\Frob_p)(f)\right),$$ where $Frob_p$ is the Frobenius element of Galois groups acting on the character space.
\end{thm}
$\hfill\Box$


\begin{thebibliography}{xxxx}
\bibitem[A]{arthur}{\sc J. Arthur}, An Introduction to the Trace Formula, in {\it Harmonic Analysis, The Trace Formula, and Shimura Varieties}, Proceedings of the Clay Mathematics Institute 2003 Summer School, The Fields Institute
Toronto, Canada, June 2–27, 2003, \textbf{Clay Mathematics Proceedings}, Vol. \textbf{4}, 2005,   1-264.

\bibitem[C]{colmez}{\sc P. Colmez}, {\it Fonctions L p-adiques}
S\'eminaire Boubaki, 51\`eme ann\'ee, 1998-99, No 851, 851-01-851-38,  Novembre 1998.

\bibitem[DQ1]{diepquynh1}{\sc D. N. Diep, D. T. P. Quynh }, {\it Automorphic representations of $\SL(2,\mathbb R)$ and quantization of fields}, {American R. J. of Math. }\textbf{1}, No 2, 2015, 25-37.

\bibitem[DQ2]{diepquynh2}{\sc D. N. Diep, D. T. P. Quynh }, {\it Poisson summation and endoscopy for $\SU(2,1)$ }, {East-West J. of Math.}, \textbf{17}, No 2,2015, 125-140.

\bibitem[DQ3]{diepquynh3}{\sc D. N. Diep, D. T. P. Quynh }, {\it Poisson summation and endoscopy for $\Sp(4,\mathbb R)$ }, {Southeast Asia Bull. of Math.}, \textbf{40}, No 5-6, 2016.

\bibitem[DQ4]{diepquynh4}{\sc D. N. Diep, D. T. P. Quynh }, {\it Poisson summation and endoscopy for $\SL(3,\mathbb R)$ }, {ArXiv:1407.6912v1 [math.QA] 13 Jul 2014}  .

\bibitem[H]{harris}{\sc M. Harris}, {\it Automorphic Galois representations and cohomology of Shimura varieties},  .

\bibitem[HLTT]{harrisetall}{\sc M. Harris, K.-W. Lan, R. Taylor, J. Thorne}, {\it On the rigid cohomology of certain Shimura varieties}


\bibitem[K]{kottwitz}{\sc R. E. Kottwitz}, Harmonic Analysis on Reductive p-adic Groups and Lie Algebras, in {\it Harmonic Analysis, The Trace Formula, and Shimura Varieties}, ibid, pp 393-522,  .

\bibitem[Ko]{kowalski}{\sc E. Kowalski}, {\it Automorphic forms, L-functions and number theory (March 12–16) Three Introductory lectures}, Preprint.


\bibitem[M]{milne}{\sc Milne}, Introduction to Shimura Varieties,  in {\it Harmonic Analysis, The Trace Formula, and Shimura Varieties}, ibid, pp 265-378.

\bibitem[W]{waldspurger}{\sc Waldspurger}, {\it },  .

\end{thebibliography}
\end{document}